\documentclass[a4paper,leqno,10pt]{amsart}
\usepackage{amsmath,amssymb,amsthm}
\usepackage[usenames,dvipsnames]{color}
\usepackage{times}
\usepackage{amsxtra}

\newtheorem{thm}{Theorem}
\newtheorem{lemma}{Lemma}
\newtheorem{propos}{Proposition}

\numberwithin{equation}{section}
\numberwithin{thm}{section}
\numberwithin{rmk}{section}
\numberwithin{lemma}{section}
\numberwithin{clr}{section}
\numberwithin{Def}{section}
\newcommand{\PP}{\mathbb{P}}
\newcommand{\EE}{\mathbb{E}}
\newcommand{\NN}{\mathbb{N}}
\newcommand{\eps}{\varepsilon}
\newcommand{\tc}{\tau_c}
\newcommand{\la}{\lambda}
\begin{document}

\title{The Deterministic and Stochastic Shallow Lake Problem}

 \author{G.~T.  Kossioris$^1$ \and M. Loulakis$^1$ \and P.~E. Souganidis$^{2}$}

 \address[G.~T.  Kossioris]{Department of Mathematics and Applied Mathematics,  Voutes Campus, University of Crete \\ 
 70013 Heraklion, Greece}
 \email{kosioris@uoc.gr}
 \thanks{$^1$This research has been co-financed by the European Union (European Social Fund ESF) and Greek national funds through the Operational Program "Education and Lifelong Learning" of the National Strategic
 Reference Framework (NSRF) 2007-2013. Research Funding Program: THALES. Investing in knowledge society through the European Social Fund, Project: Optimal Management of Dynamical Systems of the Economy and the Environment, 
 grant number MIS 377289}

 \address[Michail Loulakis]{School of Applied Mathematical and Physical Sciences \\ National Technical University of Athens \\
 15780 Athens, Greece AND Institute of Applied and Computational Mathematics- Foundation of Research and Technology Hellas\\ 700 13 Heraklion, Greece}
\email{loulakis@math.ntua.gr}

\address[P.~E. Souganidis]{Department of Mathematics,
University of Chicago \\
5734 S. University Avenue, 
Chicago, IL, 60637}
\email{souganidis@math.uchicago.edu}
\thanks{$^{2}$Partially supported by the National Science Foundation
grants DMS-1266383 and DMS-1600129 and the Office for Naval Research 
Grant
N00014-17-1-2095. }

\begin{abstract}
 We study   the welfare function of the  deterministic and stochastic shallow lake problem. We show that the welfare function is the viscosity solution of the associated Bellman equation, we establish several properties including its asymptotic
 behaviour at infinity and we present a convergent monotone numerical scheme.\\
\\
{\em AMS 2010 Mathematics Subject Classification}: 93E20, 60H10, 49L25\\
\\
{\em Keywords}: Shallow Lake, Viscosity solution, Optimal stochastic control

\end{abstract}

\maketitle
\section{Introduction}
The scope  of this work is the theoretical  study  of the welfare function
describing the economics of shallow lakes. Pollution of shallow lakes is a quite often observed phenomenon 
 because of  heavy use of fertilizers on
surrounding land and an increased inflow of waste water from human
settlements and industries.  The shallow lake system provides conflicting services as a resource, due   to the provision of ecological services of a clear lake, 
and as a waste sink, due  to agricultural activities. 
\vskip.05in

The economic analysis of the problem requires the study of an optimal control problem or a differential game  in case of  common 
property resources  by various communities; see, for example,  \cite{carpenter}, \cite{brock-starrett} and \cite{xepapadeas}. 
\vskip.05in

Typically the model is described in terms of the amount, $x(t)$, of phosphorus in algae which is assumed to evolve according to the stochastic differential equation (sde for short)
\begin{equation} \label{sldyn}
\begin{cases} dx(t)=\left( u(t)-b x(t)+\dfrac{x^{2}(t)}{x^{2}(t)+1}\right) dt+\sigma x(t)dW(t), & \\
x(0) = x. & \end{cases}
\end{equation}

 The first   term, $u(t)$,  in the drift part of the dynamics represents the exterior load of phosphorus imposed by the agricultural community, which is assumed to be positive. The second term is
  the rate of loss $bx(t)$, which consists   of sedimentation, outflow and
sequestration in other biomass.  The third term  is the rate of recycling of  phosphorus due to  
sediment resuspension resulting, for example, from waves, or oxygen depletion.  This 
rate is typically taken to be a sigmoid function; see  \cite{carpenter}. The model also assumes an uncertainty in the  recycling rate  driven by a linear multiplicative white noise with  diffusion strenght $\sigma$.   
\vskip.05in

The lake has value as a waste
sink for agriculture modeled by $\ln u$ and it provides ecological
services that decrease with the amount of phosphorus $-cx^{2}$. 
The positive parameter $c$ reflects the relative weight of this welfare
component; large $c$ gives more weight to the ecological services
of the lake.
\vskip.05in

Assuming an infinite horizon at a discount rate $\rho >0$, the total benefit is
\begin{equation}\label{Jxu}
J(x;u)=\EE\left[\int_0^\infty e^{-\rho t}\big(\ln u(t)-cx^2(t)\big)\, dt\right],
\end{equation}
where $x(\cdot)$ is the solution to (\ref{sldyn}), for a given  $u(\cdot)$, and $x(0)=x$.
Optimal management  requires to maximize the total benefit, over all  exterior loads that act as controls by the social planner. Thus  the  welfare function is
\begin{equation}\label{ihvf}
V(x) = \underset{u\in \mathfrak{U}_x} {\sup}\ J(x;u),
\end{equation}
where   $\mathfrak{U}_x$ is the set of admissible controls $u: [0,\infty) \rightarrow  \mathbb{R}^+ $  which are specified in the next section.
\vskip.05in

Dynamic programming arguments lead, under the assumption that the welfare function is decreasing,  to the  Hamilton-Jacobi-Bellman equation
$$\rho V=\left( \dfrac{x^{2}}{x^{2}+1}-bx\right)V_x-\left( \ln(-V_x)+x^{2}+1\right)+\dfrac{1}{2}\sigma^{2}x^{2}V_{xx}, \leqno\mathrm{(OHJB)}$$
and of the aims of this paper is to provide a rigorous justification of this fact.
\vskip.05in

In the deterministic case $\sigma=0$, the optimal dynamics of the problem were fully investigated by analysing the possible equilibria of the  dynamics given by Pontryagin maximum principle; see, for example, \cite{xepapadeas} and
\cite{Wagener}. The possibility to  steer the combined economic-ecological system
towards the trajectory of optimal management via optimal state-dependent taxes was also considered; see  \cite{KPXAdeZ}. 
\vskip.05in

On the other hand, there is not much in the literature about the stochastic problem. Formal  asymptotics expansions
of the solution for small  $\sigma $ for Hamilton-Jacobi-Bellman equations like  (OHJB) have been presented in \cite{GKW}. In the same paper, the  authors also give a formal  phenomenological bifurcation analysis  based on a geometric invariant quantity, along with some  numerical computations of the stochastic bifurcations
based on  (formal) asymptotics  for small $\sigma$.
\vskip.05in

The connection between  stochastic  control problems  and Hamilton-Jacobi-Bellman equation, which is based on the dynamic programing principle, has been studied extensively. 
The correct mathematical framework 
is that of the  Crandall-Lions viscosity solutions introduced in \cite{lions}; see  the review article \cite{CIL}).   The deterministic case  leads to the study of a first order Hamilton-Jacobi equation; see, for example,   
 \cite{BCD}, \cite{barles} and references therein. For the general  stochastic 
optimal control problem we refer to \cite{Kry}, \cite{PLL1} - \cite{PLL3}, \cite{FS}, \cite{JYXZ} and references therein.
\vskip.05in

The stochastic shallow lake problem has some nonstandard features and, hence, it requires some special analysis.  At first, the problem is formulated as a state constraint one on a
semi-infinite domain. Viscosity solutions with state
constraint boundary conditions were  introduced for first order equations by \cite{So} and \cite{CDL}.  For second order equations one should consult \cite{Ka},
\cite{LL} and \cite{ALL}. Apart from the left boundary condition, the correct asymptotic decay of the solution at infinity is necessary to establish a comparison result; see, for example,   \cite{IL} and
\cite{DaLL}.
\vskip.05in

The unboundedness of the controls  along with  the logarithmic term in the cost functional lead to a logarithm of the gradient variable in (OHJB).  An 
 a priori knowledge of the solution is required to guarantee that the Hamiltonian is well defined. Moreover, due the  presence of the logarithmic term it is necessary to construct an appropriate test function 
 to establish a comparison proof. The commonly used  polynomial functions, see, for example,  \cite{IL} and   \cite{DaLL}, are not useful here since they  do not yield a supersolution 
of the equation. 
\vskip.05in

In the present work we study the stochastic shallow lake problem   for a fixed $\sigma$. We first prove the necessary stochastic analysis estimates for the welfare function (\ref{ihvf}). We obtain directly various crucial
properties for the welfare function, that is,  boundary behaviour, local regularity, monotonicity and
asymptotic estimates at infinity.
\vskip.05in

We prove, including the deterministic case, that (\ref{ihvf}) is the unique decreasing constrained viscosity solution 
to (OHJB) with  quadratic growth  at infinity. The comparison theorem is  proved by considering a linearized equation and constructing a proper supersolution. Exploiting the well-posedness of the problem within the framework of constrained viscosity 
solutions we investigate the exact asymptotic behavior of the solutions at infinity. The latter is used to construct a monotone convergent numerical scheme that
along with the optimal dynamics equation can be used  to  reconstruct numerically the stochastic optimal dynamics.

%

\section{The general setting and the main results}

We assume that there exists a filtered probability space $( \Omega, \mathcal{F},\{ \mathcal{F}_t\}_{t\ge 0},\PP )$
satisfying the usual conditions, and a Brownian motion $W(\cdot)$ defined on that space. 
An admissible control   $u(\cdot)\in\mathfrak{U}_x$ is an $\mathcal{F}_t$-adapted, $\PP$-a.s. locally integrable process with values in $U=(0,\infty)$, satisfying 
\begin{equation}\label{ac_constraint}
 \mathbb{E} \left[ \int_{0}^{\infty}e^{-\rho t}\ln u(t)dt \right] < \infty,
\end{equation}
such that  the problem (\ref{sldyn}) has a unique strong solution $x(\cdot)$.
\vskip.05in

 The shallow lake problem has an infinite horizon. Standard arguments  based on the dynamic programming principle   (see \cite{FS}, Section III.7)  suggest that the welfare function  $V$  given by (\ref{ihvf}) satisfies  the  HJB equation
\begin{equation}\label{HJB}
\rho V=\sup_{u\in U } G(x,u,V_x,V_{xx}),
\end{equation}
with  $G $  defined by
\begin{equation}\label{hG}
G(x,u,p,P)=\dfrac{1}{2}\sigma^{2}x^{2}P+\left( u-bx+\dfrac{x^{2}}{x^{2}+1}\right)p +\ln u-x^{2}.
\end{equation}

One difficulty in  the study of this problem is related to the fact that the control functions $u$ take values in the unbounded set $ \mathbb{R}^+$ so that supremum in (\ref{HJB}) 
might take infinite values. Indeed, when $U=\mathbb{R}^+$, setting 
$$ H(x,p,P)= \sup_{u\in \mathbb{R}^+ } G(x,u,p,P),$$
we find 
\begin{equation}\label{HJB-unbounded}
H(x,p,P)=
  \left\{
  \begin{array}{lr}  \left( \dfrac{x^{2}}{x^{2}+1}-bx\right)p-\left( \ln(-p)+x^{2}+1\right)+\dfrac{1}{2}\sigma^{2}x^{2}P & \,\, if \,\, p<0, \\
                                                     + \infty & \,\, if \,\, p \geq 0. \end{array} \right.
 \end{equation}

 It is natural to expect that  since shallow lake looses its value  with a higher  concentration  of phosphorus,  the welfare function is decreasing
 as the initial state of phosphorus increases. Assuming that $V_x<0$,  (\ref{HJB}) becomes (OHJB).
\vskip.05in 

Since the problem is set in  $(0, \infty)$, it is necessary to introduce boundary conditions   guaranteeing the well-posedness of  the corresponding boundary value problem. 
\vskip.05in

Given the possible  degeneracies  of Hamilton-Jacobi-Bellman equations at $x=0$, the right framework 
is that of continuous viscosity solutions in which boundary conditions are considered in the weak sense.
Since at the boundary point $x=0$
\begin{equation}\label{control}
\underset{u\in U} {\inf} \left\{ - u+b x-\dfrac{x^{2}}{x^{2}+1} \right\} < 0,
\end{equation}
that is,  there always exists a control that
can  drive the system inside $(0,\infty)$, the problem should  be considered as a state constraint one on the interval $[0,\infty)$, meaning that $V$ is a subsolution in $[0,\infty)$ and a supersolution in $(0,\infty)$. 
\vskip.05in

Next we present the main results of the paper. The proofs are given in Section 4. The first is about the relationship between the welfare function and (OHJB). 

\begin{thm}\label{constrained_vs}
If $\sigma^2 < \rho +2b$, the welfare function $V$  is a continuous in $[0,\infty) $ constrained viscosity solution of the  equation {\rm (OHJB)} in  $[0,\infty)$. 
\end{thm}
\noindent
The  second  result is the following  comparison principle for solutions of (OHJB).   
\begin{thm}\label{comparison} Assume that  $u\in C([0,\infty))$ is a bounded from above  strictly decreasing subsolution of ${\rm (OHJB)}$ in $[0, \infty)$ and $v\in C([0,\infty))$ is a bounded from above strictly decreasing supersolution of
${\rm (OHJB)}$ in $(0, \infty)$ such that  $v \geq -c(1+x^2)$ and  $Du\leq -\frac{1}{c^*}$
in the viscosity sense,
for $c$, $c^*$  positive constants. Then
$u \leq v$ in $[0, \infty)$. 
\end{thm}
\noindent

The next theorem describes the exact asymptotic behavior of (\ref{ihvf}) at $+\infty$. Let 
\begin{equation*}\label{K}
A=\frac{1}{\rho+2b-\sigma^2}\,\,\, \mbox{ and } \,\,\, K=\frac{1}{\rho}\left(\frac{2b+\sigma^2}{2\rho}-\frac{A(\rho+2b)}{(b+\rho)^2} -1\right).
\end{equation*}

\begin{thm}\label{asymptotic}
As $x\rightarrow \infty$,
\begin{equation}\label{asympto_behav}
V(x)= -A\left(x+\frac{1}{b+\rho}\right)^2-\frac{1}{\rho}\ln \left[2A(x+\frac{1}{b+\rho})\right]+K+o(1). 
\end{equation}
\end{thm}
\vskip.05in

An important ingredient of the analysis is the following proposition which collects some key properties of $V$ that  are used in 
the proofs of Theorem \ref{constrained_vs} and Theorem \ref{asymptotic} and  show that $V$ satisfies the assumptions of Theorem \ref{comparison}. The proof is presented in Section 3. 
\begin{propos}\label{vprop}
Suppose $\sigma^2<\rho+2b$.\\[1mm]
(i) There exist constants $K_1, K_2>0$, such that, for any $x\ge 0$, we have
\begin{equation}\label{Vbounds}
K_1\ \le\ V(x)+A\left(x+\frac{1}{b+\rho}\right)^2+\frac{1}{\rho}\ln \left(x+\frac{1}{b+\rho}\right)\ \le\ K_2.
\end{equation}
\noindent
(ii) There exist  $C>0$ and  $\Phi:[0,+\infty)\rightarrow \mathbb{R} $ increasing such that, for any $x_1,x_2\in [0,+\infty)$ with $x_1<x_2$,
\begin{equation}\label{DVbounds}
-\Phi(x_2)\le \frac{V(x_2)-V(x_1)}{x_2-x_1}\le -C.
\end{equation}
(iii) $V$ is differentiable at zero and 
\begin{equation}
\ln\big(-V'(0)\big)+\rho V(0)+1=0.
\label{bc}
\end{equation} 
\end{propos}

It is shown in the next section that the assumption $\sigma^2<\rho+2b$ is necessary, otherwise $V(x)=-\infty.$  
\vskip.05in

Relation (\ref{asympto_behav}) is important for the numerical approximation of (\ref{ihvf}). Given that the  computational domain   is finite, the correct asymptotic behavior
of (\ref{ihvf}) at $+\infty$ is necessary for an accurate computation of $V$.
In this connection, in Section 5 we present a monotone numerical scheme approximating (\ref{ihvf}). 

 \section{The proof of Proposition \ref{vprop}}
\noindent
 {\bf Properties of the dynamics.} The first result states that, if $x\geq0$,  the solution to (\ref{sldyn}) stays nonnegative.
\begin{lemma}\label{positive}
If $x\ge 0$, $u\in\mathfrak{U}_{x}$, and $x(\cdot)$ is the solution to (\ref{sldyn}), then
$
\PP\big[x(t)\ge 0,\ \forall t\ge 0\big]=1.
$
\end{lemma}
\noindent
{\em Proof:}\ Elementary stochastic analysis calculations yield that 
\begin{equation}
x(t)=xZ_t+\int_0^t \frac{Z_t}{Z_s}\left(u(s)+\frac{x^2(s)}{1+x^2(s)}\right) ds,
\label{implicit}
\end{equation}
where 
\begin{equation}\label{fundsol}
Z_t=e^{\sigma W_t-(b+\frac{\sigma^2}{2})t},
\end{equation}
and the claim is now obvious, since $u$ takes positive values.\hfill$\Box$\\[1mm]

The next assertion is  that the set of admissible controls is actually independent of the starting point $x$. 
\begin{lemma}\label{uxisu} 
  For all $x,y\ge 0$,  $\mathfrak{U}_x=\mathfrak{U}_y.$
\end{lemma}
\noindent
{\em Proof:} Fix  $u\in\mathfrak{U}_x$ and $x\in[0,\infty)$ and let $x(\cdot)$  the unique strong solution to (\ref{sldyn}) with $x(0)=x$, and, for any $y\ge 0$, consider the sde 
\begin{equation}\label{U_x-indep}
\begin{cases} dw(t)=\Big\{-bw(t)-\frac{(x(t)-w(t))^2}{1+(x(t)-w(t))^2}+\frac{x(t)^2}{1+x(t)^2}\Big\}dt+\sigma w(t)dW(t),&\\
w(0)=x-y,\end{cases}
\end{equation}
and note that the coefficients  are Lipschitz and grow at most linearly.
It follows that (\ref{U_x-indep}) has a unique strong solution defined for all $t\ge 0$.  It is easy to see now that the process $y(t)=x(t)-w(t)$ satisfies (\ref{sldyn}) with $y(0)=y$. Moreover, the uniqueness of $y(\cdot)$ follows from that of $x(\cdot)$, so $u\in\mathfrak{U}_y$. \hfill$\Box$\\[2mm]
In view of Lemma \ref{uxisu} we will denote the set of admissible controls by $\mathfrak{U}$, regardless of the starting point $x$  in (\ref{sldyn}).

\begin{lemma}\label{monotone}
Suppose $x(\cdot),\ y(\cdot)$ satisfy (\ref{sldyn}) with controls $u_1,\ u_2\in\mathfrak{U}$, respectively, and $x(0)=x$, $y(0)=y$.
If $x\le y$ and $\PP\big[u_1(t)\le u_2(t),\ \forall t\ge 0\big]=1$, then
\[
\PP\big[y(t)-x(t)\ge (y-x)Z_t,\ \forall t\ge 0\big]=1.
\]
\end{lemma}
\noindent
{\em Proof:} The proof is an immediate consequence of (\ref{implicit}) and Gronwall's inequality, since  $w(t)=y(t)-x(t)$ satisfies
\[
w(t)=(y-x)Z_t+\int_0^t \frac{Z_t}{Z_s}\big(u_2(s)-u_1(s)+w(s)\frac{y(s)+x(s)}{\big(1+y^2(s)\big)\big(1+x^2(s)\big)}\big)\ ds.
\]
 $\hfill\Box$\\[5mm]
\noindent
{\bf Properties of the welfare function.} The results here follow from the properties of (\ref{sldyn}).  In the rest of the paper, for the shake of convenience, we  
assume that $c=1$. Throughout this section, we  refer to quantities depending only on $\rho,b$ and $\sigma^2$  as constants.
\vskip.05in

We remark that, if $\sigma^2\ge \rho+2b$, then $V(x)\equiv-\infty$. Indeed, when $u\in\mathfrak{U}$ and $x(\cdot)$ satisfies (\ref{sldyn}), 
(\ref{implicit}) implies that $\PP\big[x(t)\ge M_t(u),\ \forall t\ge 0\big]=1$, with 
\begin{equation}\label{Mf}
M_t(u)=\int_0^t\frac{Z_t}{Z_s}u(s)\ ds,
\end{equation}
where $\{Z_t\}_{t\ge 0}$ is as in \eqref{fundsol}. 

\vskip.05in 
In view of this observation we will hereafter assume that 
\[
 \sigma^2<\rho+2b.
\]

We will first prove three lemmata before we proceed with the proof of Proposition \ref{vprop}.
\begin{lemma}\label{Vdec}
 The function $x\mapsto V(x)+Ax^2$ is decreasing on $[0,+\infty)$. 
\end{lemma}

\noindent
{\em Proof:} Fix $x_1,x_2\ge 0$ with $x_1\le x_2$. It suffices to show that, for any  control $u\in\mathfrak{U}$,
\[
J(x_2;u)+Ax_2^2\le J(x_1;u)+Ax_1^2.
\]

Since this holds trivially if $J(x_2;u)=-\infty$, we may assume that $J(x_2;u)>-\infty$. 
\vskip.05in

Consider now the solutions $x_1(\cdot),x_2(\cdot)$ to (\ref{sldyn}) 
with initial conditions $x_1, x_2$ and a common control $u\in\mathfrak{U}$. Lemma \ref{monotone} implies that, $\PP$-a.s.  and for all $t\geq0$ and $i=1,2,$
\[
x_i(t)\ge x_iZ_t, \qquad\text{and}\qquad x_2(t)-x_1(t)\ge (x_2-x_1)Z_t.
\]

Note that since $u\in\mathfrak{U}$ and $J(x_2;u)>-\infty$,
\[
\int_0^\infty e^{-\rho t}x_2^2(t)\,dt<+\infty\Rightarrow\int_0^\infty e^{-\rho t}x_1^2(t)\,dt<+\infty\Rightarrow J(x_1;u)>-\infty.
\]

In particular,
\begin{align*}
J(x_2;u)\!-\!J(x_1;u)&=\EE\left[ \int_{0}^{\infty}\!\!\!\!e^{-\rho t}\big(\ln u(t)-x_2^{2}(t)\big)dt -\!\! \int_{0}^{\infty}\!\!\!\!e^{-\rho t}\big(\ln u(t)-x_1^{2}(t)\big)dt\right]\\
&=-\EE\left[ \int_{0}^{\infty} e^{-\rho t}\big(x_2(t)-x_1(t)\big)\big(x_2(t)+x_1(t)\big)dt \right]\\
&\le -(x_2^2-x_1^2)\int_0^\infty e^{-\rho t} \EE\big[Z_t^2\big]\ dt= -A(x_2^2-x_1^2).
\end{align*}
\hfill$\Box$
\begin{lemma}\label{V0fin}
The welfare function at zero satisfies
$
V(0) \le\frac{1}{\rho}\ln\left(\frac{b+\rho}{\sqrt{2e}}\right).
$
\end{lemma}
\noindent
{\em Proof:} Recall that, for any $u\in\mathfrak{U}$,  $x(t)\ge M_t(u)$.
Using Jensen's  inequality and part (i) of Lemma \ref{Mest}, we find
\begin{align*}
\mathbb{E}\left[\int_0^\infty e^{-\rho t}\ln u(t)\ dt\right] &\le\frac{1}{\rho}\ln\mathbb{E}\left[\int_0^\infty \rho e^{-\rho t}u(t)\ dt\right] \\
&=\frac{1}{\rho}\ln\mathbb{E}\left[\int_0^\infty \rho(\rho+b) e^{-\rho t}M_t(u)\ dt\right] \\
&\le \frac{\ln(b+\rho)}{\rho}+\frac{1}{\rho}\ln\mathbb{E}\left[\int_0^\infty \rho e^{-\rho t}x(t)\ dt\right] \\
&\le \frac{\ln(b+\rho)}{\rho}+\frac{1}{2\rho}\ln\mathbb{E}\left[\int_0^\infty \rho e^{-\rho t}x^2(t)\ dt\right].
\end{align*}

In view of (\ref{ac_constraint}), we need only consider  $u\in\mathfrak{U}$ such that $D=\mathbb{E}\left[\int_0^\infty e^{-\rho t}x^2(t)\ dt\right]<\infty$.
Then
\[
\mathbb{E} \left\{ \int_{0}^{\infty}e^{-\rho t}\big[\ln u(t)-x^{2}(t)\big]dt \right\}\le  \frac{\ln(b+\rho)}{\rho}+\frac{\ln(\rho D)}{2\rho}-D\le\frac{1}{\rho}\ln\left(\frac{b+\rho}{\sqrt{2e}}\right),
\]
and the assertion holds.\hfill$\Box$\\

It follows from Lemma \ref{Vdec} and Lemma \ref{V0fin} that $V<+\infty$ in $[0,\infty)$.
\noindent
The next result is a special case of the dynamic programming principle. 
\begin{lemma}\label{specialdpp}
Fix $x_1,x_2\in[0,\infty)$ with $x_1<x_2$, and, for  $u\in\mathfrak{U}$, let $x(\cdot)$ be the solution to (\ref{sldyn}) with control $u$ and $x(0)=x_1$.  If $\tau_u$ is 
the hitting time of $x(\cdot)$ on $[x_2,+\infty)$, that is,
\[
\tau_u=\inf\{t\ge 0: x(t)\ge x_2\},
\]
then

\begin{equation}\label{babydpp}
V(x_1)=\sup_{u\in\mathfrak{U}}\EE\left[\int_0^{\tau_u} e^{-\rho t}\big(\ln u(t)-x^2(t)\big)\, dt+e^{-\rho\tau_u}V(x_2)\right].
\end{equation}
\end{lemma}
\noindent
{\em Proof:} We have
\begin{align*}
J(x_1;u)&=\EE\left[\int_0^{\tau_u}\hspace{-2mm}e^{-\rho t}\big(\ln u(t)-x^2(t)\big)\, dt\right]
+\EE\left[\int_{\tau_u}^\infty \hspace{-2mm}e^{-\rho t}\big(\ln u(t) -x^{2}(t)\big)\,dt;\, \tau_u<+\infty\right].
\end{align*}

Conditioning on the $\sigma$-field $\mathcal{F}_{\tau_u}$, and applying the strong Markov property, the rightmost term becomes
\begin{align*}
\EE\left[e^{-\rho\tau_u}\EE\left[\int_{\tau_u}^\infty \hspace{-2mm}e^{-\rho (t-\tau_u)}\big(\ln u(t) -x^{2}(t)\big)\,dt\ \Big| \mathcal{F}_{\tau_u}\right];\, \tau_u<+\infty\right]\le \EE\left[e^{-\rho\tau_u}\right]V(x_2),
\end{align*}
since on the event $\{\tau_u<+\infty\}$, $x(\tau_u+\cdot)$ satisfies (\ref{sldyn}) with initial condition $x(\tau_u)=x_2$ and control $u(\tau_u+\cdot)$. Taking the supremum over $u\in\mathfrak{U}$ we
see that the left hand side of (\ref{babydpp}) is less than or equal the right hand one.
\vskip.05in
For the reverse inequality, take any $u\in\mathfrak{U}$ and consider (\ref{sldyn}) driven by the Brownian motion 
$B(t)=W(\tau_u+t)-W(\tau_u)$, and,  for $\eps>0$, choose  a control $u_\eps$ such that 
\[
V(x_2)<J(x_2;u_\eps)+\eps.
\]

Define now the new control $u_*\in\mathfrak{U}$ as 
\[
u_*(t)=\begin{cases} u(t) & \text{ for } t\le\tau_{u}\\ u_\eps(t-\tau_u) & \text{ for } t> \tau_u.\end{cases}
\]

Just as in the proof of the upper bound we get
\begin{align*}
V(x_1)&\ge J(x_1;u_*)=\EE\left[\int_0^{\tau_u} e^{-\rho t}\big(\ln u(t)-x^2(t)\big)\, dt+e^{-\rho\tau_u}J(x_2;u_\eps)\right]\\
&>\EE\left[\int_0^{\tau_u} e^{-\rho t}\big(\ln u(t)-x^2(t)\big)\, dt+e^{-\rho\tau_u}V(x_2)\right]-\eps,
\end{align*}
which concludes the proof.\hfill$\Box$\\[2mm]

We next give the proof of Proposition \ref{vprop}, which is subdivided in several parts. \\[1mm]

\noindent
{\bf Proof of Proposition \ref{vprop}:} 
\noindent
{\em Proof  lower bound, part (i):}\quad The claim will follow by choosing the  control $u(t)=\frac{1+x(t)}{1+x^2(t)}$, which is clearly admissible. 
Then, (\ref{implicit}) gives
\begin{equation}
x(t)=xZ_t+\int_0^t \frac{Z_t}{Z_s}\ \Big(1+\frac{x(s)}{1+x^2(s)}\Big)ds=xZ_t+M_t(1)+\int_0^t \frac{Z_t}{Z_s}\ \frac{x(s)}{1+x^2(s)}ds
\label{testX}
\end{equation}
and, hence,
\begin{align*}
x^2(t)&=x^2Z_t^2+2xZ_tM_t(1)+M_t^2(1)+\Big(\int_0^t \frac{Z_t}{Z_s}\ \frac{x(s)}{1+x^2(s)}ds\Big)^2\\
&\quad+M_t(1)\int_0^t \frac{Z_t}{Z_s}\ \frac{2x(s)}{1+x^2(s)}ds+\int_0^t \frac{Z_t^2}{Z_s}\ \frac{2x x(s)}{1+x^2(s)}ds.
\end{align*}

To estimate the rightmost term from above, note that, in view of (\ref{testX}),  $x\le x(s)Z_s^{-1}$,  while for the third and fourth terms of the sum we use that
$\frac{x(s)}{1+x^2(s)}\le \frac{1}{2}$. It follows that
\begin{equation}
x^2(t)\le x^2Z_t^2+2xZ_tM_t(1)+\frac{9}{4}M_t^2(1)+2\int_0^t \frac{Z_t^2}{Z_s^2}\ ds,
\label{testX2}
\end{equation}

It is easy to see that
\begin{equation}
\EE\left[\int_0^\infty e^{-\rho t}x^2Z_t^2 dt\right]=x^2\int_0^\infty e^{-(\rho+2b-\sigma^2)t}dt=Ax^2,
\label{Xsq}
\end{equation}

Lemma \ref{Mest} (ii) gives
\begin{align}\label{X1}
\EE\left[\int_0^\infty e^{-\rho t} 2xZ_t M_t(1) dt\right]&=2Ax\int_0^\infty e^{-\rho t}\EE\big[Z_t\big]dt= \frac{2Ax}{\rho+b},
\end{align}
while Lemma \ref{Mest} (i) and Lemma \ref{Mest} (iii)  yield
\begin{equation}
\int_0^\infty e^{-\rho t}\EE\big[M_t^2(1)\big]\ dt = 2A\int_0^\infty e^{-\rho t}\EE\big[M_t(1)\big]\ dt =\frac{2A}{\rho(\rho+b)}.
\label{Mtsq}
\end{equation}

We also have
\begin{align}\label{lastterm}
\int_0^\infty \hspace{-2mm}e^{-\rho t}\int_0^t\EE\left[\frac{Z_t^2}{Z_s^2}\right]ds\ dt&=\int_0^\infty \hspace{-2mm}e^{-\rho t}\int_0^t e^{(\sigma^2-2b)(t-s)}ds\ dt=\frac{A}{\rho}.
\end{align}

Using the last four observations in  (\ref{testX2}), we find  for some constant $B$, 
\begin{equation}
\int_0^\infty e^{-\rho t}\EE\big[x^2(t)\big]\ dt \le A\left[\big(x+\frac{1}{\rho+b}\big)^2+B\right].
\label{Xener}
\end{equation}

 On the other hand, using that, for all $x\ge 0$, $(1+x)^2\ge (1+x^2)$,  and Jensen's inequality,  we find
\begin{align*}
\int_0^\infty e^{-\rho t}\EE\big[\ln u(t)\big]\ dt &\ge-\int_0^\infty e^{-\rho t}\EE\big[\ln\big(1+x(t)\big)\big]\ dt \\
&\ge -\frac{1}{\rho}\ln\left(\int_0^\infty \rho e^{-\rho t} \Big(1+\EE\big[x(t)\big]\Big) dt\right)\\
&=-\frac{1}{\rho}\ln\left(1+\rho\int_0^\infty  e^{-\rho t} \EE\big[x(t)\big] dt\right).
\end{align*}

By  (\ref{testX}) it follows that
$\displaystyle
\EE\big[x(t)\big]\le x\EE\big[Z_t\big]+\frac{3}{2}\EE\big[M_t(1)\big]=xe^{-bt}+\frac{3}{2}\EE\big[M_t(1)\big].
$

Hence, using Lemma \ref{Mest} (i), we obtain
\[
\int_0^\infty e^{-\rho t}\EE\big[\ln u(t)\big]\ dt \ge -\frac{1}{\rho}\ln\left(1+\frac{\rho x}{\rho+b}+\frac{3}{2(\rho+b)}\right).
\]

The preceding estimate and (\ref{Xener}) together imply that, for some suitable constant $K_1$,
\begin{align*}
V(x)&\ge J(x;u) =\EE\big[\int_0^\infty e^{-\rho t}\big(\ln u(t)-x^2(t)\big)\ dt \big]\\
&\ge -A\left(x+\frac{1}{b+\rho}\right)^2-\frac{1}{\rho}\ln \left(x+\frac{1}{b+\rho}\right)+K_1.
\end{align*}

\noindent
{\em Proof of the upper bound, part (i) :}\quad In view of Lemma \ref{Vdec} and Lemma \ref{V0fin}, it suffices to find  $K_2>0$, such that the asserted inequality holds for $x\ge 1$.
\vskip.05in

Fix $u\in\mathfrak{U}$. Then 
\begin{align*}
x^2(t)&\ge x^2Z_t^2+2xZ_t^2\int_0^t\frac{1}{Z_s}\left(u(s)+\frac{x^2(s)}{1+x^2(s)}\right) ds\\
&=x^2Z_t^2+2xZ_tM_t(1+u)-\int_0^t\frac{Z_t^2}{Z_s^2}\frac{2xZ_s}{1+x^2(s)} ds.
\end{align*}

Since
\[
1+x^2(t)\ge 1+x^2Z_t^2\ge 2xZ_t,
\]
so we can further estimate $x^2(t)$ from below by
\begin{align}\label{mtuest}
x^2(t)&\ge x^2Z_t^2+2xZ_tM_t(1)+2xZ_tM_t(u)-\int_0^t \frac{Z_t^2}{Z_s^2} ds.
\end{align}

Using the elementary inequality that $\ln a\le ab-\ln b-1$, which holds for all $a,b>0$, and Lemma \ref{Mest} (ii), we obtain, for some B,
\begin{align*}
&\int_0^\infty e^{-\rho t}\EE\big[\ln u(t)\big]dt \le \EE\left[\int_0^\infty e^{-\rho t} \Big\{2AxZ_tu(t)-\ln\big(2AxZ_t\big)\Big\}\ dt\right]\nonumber\\
&\qquad=\EE\left[\int_0^\infty\hspace{-2mm} e^{-\rho t} 2xZ_tM_t(u)\, dt\right]-\frac{\ln(2Ax) }{\rho}+\frac{2b+\sigma^2}{2\rho^2}\nonumber\\
&\qquad\le \EE\left[\int_0^\infty \hspace{-2mm}e^{-\rho t} \Big(x^2(t)-x^2Z_t^2-2xZ_tM_t(1)+\int_0^t\frac{Z_t^2}{Z_s^2}\, ds\Big)\, dt\right]-\frac{\ln x + B}{\rho}\nonumber,
\end{align*}
where in the final step we have used (\ref{mtuest}). 
\vskip.05in

In view of (\ref{Xsq}), (\ref{X1}) and (\ref{lastterm}), for every $u\in\mathfrak{U}$ there exists  $K_2>0$ such that
\[
J(x;u)\le -A\left(x+\frac{1}{b+\rho}\right)^2-\frac{1}{\rho}\ln \big(x+\frac{1}{\rho+b}\big)+K_2.
\]

The assertion now follows by taking the supremum over $u\in\mathfrak{U}$.\\[2mm]
\noindent
{\em Proof of the lower bound, part (ii):} Fix $x_1,x_2$ as in the statement. It follows from  Lemma \ref{specialdpp} that for any $\epsilon>0$ there exists a control $u_\eps\in\mathfrak{U}$ such that 
\begin{equation}\label{dppe}
V(x_1)\le \EE\left[\int_0^{\tau_{\eps}}e^{-\rho t}\ln u_\eps(t)\ dt\right]+\EE\big[e^{-\rho\tau_{\eps}}\big]V(x_2)+\eps(x_2^2-x_1^2),
\end{equation}
where $\tau_\eps$ is the hitting time on $[x_2,+\infty)$ of the solution $x_\eps(\cdot)$ to (\ref{sldyn}) with $x(0)=x_1$ and control $u_\eps$.
\vskip.05in

Using the elementary inequality
\[
\ln u_\eps(t)\le \ln c+\frac{u_\eps(t)}{c}-1, \quad\text{with } c=e^{\rho V(x_2)+1},
\]
we find
\begin{align}\label{DVcont}
V(x_1)-V(x_2)&\le \frac{1}{c}\EE\left[\int_0^{\tau_\eps} e^{-\rho t}u_\eps(t)\ dt\right]+\eps(x_2^2-x_1^2).
\end{align}

To conclude it suffices to show that
\begin{equation}
\frac{1}{c}\EE\left[\int_0^{\tau_\eps}e^{-\rho t}u_\eps(t)\ dt\right]\le \Phi(x_2)(x_2-x_1).
\label{suffu}
\end{equation}

To this end, we apply It\^{o}'s rule to the semimartingale $Y_t=e^{-\rho t+\gamma x_\eps(t)}$, where $\gamma>0$ is a constant to be determined, and find
\begin{align*}
Y_t-e^{\gamma x_1}&=\int_0^t Y_s \big(-\rho\, ds+ \gamma\, dx_\eps(s)+\frac{\gamma^2}{2}\, d\langle x_\eps\rangle_s\big)\\
&=\int_0^t Y_s\Big(-\rho+\gamma\big(u_\eps(s)-bx_\eps(s)+\frac{x_\eps^2(s)}{1+x_\eps^2(s)}\big)+\frac{\gamma^2\sigma^2 x_\eps^2(s)}{2}\Big)ds+M_t,
\end{align*}
where $M_t$ stands for the martingale $\gamma\sigma \int_0^t x_\eps(s)\, dW(s)$. 
\vskip.05in
Next we  apply the optional stopping theorem for the bounded stopping time $\tau_N=\min\{\tau_\eps,N\}$, with $N\in\NN$, to find 
\begin{align*}
\EE\big[Y_{\tau_N}\big]-e^{\gamma x_1}&=\EE\left[\int_0^{\tau_N} \hspace{-2mm}Y_s\Big(-\rho+\gamma\big(u_\eps(s)-bx_\eps(s)+\frac{x_\eps^2(s)}{1+x_\eps^2(s)}\big)+\frac{\gamma^2\sigma^2 x_\eps^2(s)}{2}\Big)ds\right].
\end{align*}

Since  $0\le x_\eps(s) \le x_2$ in $[0,\tau_\eps]$, 
\begin{align*}
e^{\gamma x_2}\EE\big[e^{-\rho\tau_N}\big]-e^{\gamma x_1}&\ge \EE\left[\int_0^{\tau_N} \hspace{-2mm}Y_t\Big(-\rho+\gamma\big(u_\eps(t)-bx_\eps(t)\big)+\frac{\gamma^2\sigma^2 x_\eps^2(t)}{2}\Big)dt\right],
\end{align*}
and
\begin{align*}
e^{\gamma x_2}-e^{\gamma x_1}&\ge \gamma\EE\left[\int_0^{\tau_N}\hspace{-2mm}e^{-\rho t}u_\eps(t)\, dt\right]+\EE\left[\int_0^{\tau_N} \hspace{-2mm}Y_t\Big(-b\gamma x_\eps(t)+\frac{\sigma^2\gamma^2 x_\eps^2(t)}{2}\Big)dt\right].
\end{align*}

Note that the term in the parenthesis above is nonnegative if $\gamma x_\eps(t) \ge 2b/\sigma^2$, and greater than or equal to $-b^2/2\sigma^2$ in any case. Hence, 
\[
e^{\gamma x_2}-e^{\gamma x_1}\ge \gamma\EE\left[\int_0^{\tau_N}\hspace{-2mm}e^{-\rho t}u_\eps(t)\, dt\right]-\frac{b^2e^{\frac{2b}{\sigma^2}}}{2\sigma^2}\EE\left[\int_0^{\tau_N} \hspace{-2mm}e^{-\rho t}dt\right].
\]

Letting $N\to\infty$  we get
\begin{equation}\label{coerce}
e^{\gamma x_2}-e^{\gamma x_1}\ge \gamma\EE\left[\int_0^{\tau_\eps}\hspace{-2mm}e^{-\rho t}u_\eps(t)\, dt\right]-\frac{b^2e^{\frac{2b}{\sigma^2}}}{2\sigma^2}\EE\left[\int_0^{\tau_\eps} \hspace{-2mm}e^{-\rho t}dt\right].
\end{equation}

With $\gamma$ still at our disposal,  to show (\ref{suffu}) it  suffices to control the term $\EE\left[\int_0^{\tau_\eps} e^{-\rho t}dt\right]$
by $\EE\left[\int_0^{\tau_\eps}e^{-\rho t}u_\eps(t)\, dt\right]$. 
\vskip.05in

Since without loss of generality we may assume that $\eps<A$,  Lemma \ref{Vdec} and (\ref{dppe}) give
\[
0\le V(x_1)-V(x_2)-\eps(x_2^2-x_1^2)\le \EE\left[\int_0^{\tau_\eps}e^{-\rho t}\ln u_\eps(t)\ dt\right]-\rho V(x_2)\EE\left[\int_0^{\tau_\eps}e^{-\rho t}\, dt\right].
\]

Jensen's inequality then implies that
\[
\EE\left[\int_0^{\tau_\eps}e^{-\rho t}u_\eps(t)\, dt\right] \ge e^{\rho V(x_2)}\ \EE\left[\int_0^{\tau_\eps}e^{-\rho t}\, dt\right],
\]
and (\ref{coerce}) gives, with $C(x_2)=\frac{b^2}{2\sigma^2}e^{\frac{2b}{\sigma^2}-\rho V(x_2)}$,
\begin{equation}\label{coerce2}
\gamma e^{\gamma x_2}(x_2-x_1)\ge e^{\gamma x_2}-e^{\gamma x_1}\ge (\gamma-C(x_2))\ \EE\left[\int_0^{\tau_\eps}e^{-\rho t}u_\eps(t)\, dt\right].
\end{equation}

Even though (\ref{suffu}), and hence the assertion of the Theorem, follows now with a suitable choice of $\gamma$, we will optimize the preceding inequality for later use. 
Choosing $\gamma=q(x_2)/x_2$ in \eqref{coerce2}, where
\begin{equation}\label{qx2}
q(x_2)=\frac{x_2C(x_2)}{2}+\sqrt{\left(\frac{x_2C(x_2)}{2}\right)^2+x_2C(x_2)},
\end{equation}
we obtain
\[
\EE\left[\int_0^{\tau_\eps}e^{-\rho t}u_\eps(t)\, dt\right] \le (x_2-x_1)\big(1+q(x_2)\big)e^{q(x_2)}.
\]

Letting $\eps \rightarrow 0$, (\ref{DVcont}) yields
\begin{equation}\label{DVlower}
\frac{V(x_2)-V(x_1)}{x_2-x_1}\ge -\big(1+q(x_2)\big)e^{q(x_2)-1-\rho V(x_2)},
\end{equation}
the claim now follows.\hfill$\Box$ \\[2mm]

\noindent
{\em Proof of the  upper bound, part (ii):}  In view of Lemma \ref{Vdec}, it suffices to assume that $x_2\le b$, since otherwise we have
\[
V(x_2)-V(x_1)\le -A(x_2^2-x_1^2)< -Ab(x_2-x_1).
\]

For a positive constant $c$, choose  a $u_c\in\mathfrak{U}$ that is constant an equal to $c$ up to time $\tau_c=\tau_{u_c}$.
Then, Lemma \ref{specialdpp} yields
\[
V(x_1)\ge \frac{\ln c-x_2^2}{\rho}(1-\EE\big[e^{-\rho\tc}\big])+\EE\big[{e^{-\rho\tc}}\big] V(x_2),
\]
or equivalently,
\begin{equation}
\big(V(x_2)-V(x_1)\big)\EE\big[e^{-\rho\tc}\big]\le -\big(\ln c-\rho V(x_1)-x_2^2)\ \EE\left[\int_0^{\tc}e^{-\rho t}\, dt\right].
\label{ub}
\end{equation}

Consider now the solution $x_c(\cdot)$ to (\ref{sldyn}) with $x(0)=x_1$ and control $u_c$. Applying It\^o's formula to $e^{-\rho t}x_c(t)$, 
followed by the optional stopping theorem for the bounded stopping time $\tau_N=\tc\wedge N$, we get
\begin{align}
\EE\big[e^{-\rho \tau_N}x_c(\tau_N)\big]-x_1&=\EE\left[\int_0^{\tau_N} e^{-\rho t}\big(c-(b+\rho)x_c(t)+\frac{x_c^2(t)}{1+x_c^2(t)}\big)\, dt\right].
\label{lastone}
\end{align}

The leftmost term of \eqref{lastone} is equal to \(x_2\EE\big[e^{-\rho \tc};\tc\le N\big]+e^{-\rho N}\EE\big[x_c(\tau_N); \tc>N\big]\). 
\vskip.05in

On the other hand, since we have assumed that $x_2\le b$, we have $x_c(t)\le b$ up to time $\tc$. Thus, the right hand side 
of \eqref{lastone} is bounded by \( \EE\left[\int_0^{\tau_N}e^{-\rho t} c\ dt\right]\). 
\vskip.05in

Letting $N\to\infty$ in $\eqref{lastone}$, by the monotone convergence theorem,  we have
\[
x_2 \EE\big[e^{-\rho \tc}\big]-x_1\le c\,\EE\big[\int_0^{\tc} e^{-\rho t}dt\big] \Longleftrightarrow    (x_2-x_1)\EE\big[e^{-\rho\tau_c}\big]\le (c+\rho x_1)\ \EE\big[\int_0^{\tc} e^{-\rho t}dt\big].
\]

Substituting this in (\ref{ub}) and choosing $\ln c=\rho V(x_1)+1+x_2^2$, we find
\begin{equation}\label{lowb}
V(x_2)-V(x_1)\le -(x_2-x_1)\left(e^{\rho V(x_1)+1+x_2^2}+\rho x_1\right)^{-1} \!\!\!.
\end{equation}

The assertion now follows letting $C=Ab\wedge \left(e^{\rho V(0)+1+b^2}+\rho b\right)^{-1}>0$. \hfill$\Box$ \\[2mm]

\noindent
{\em Proof of part (iii):} It follows from (\ref{lowb}) that, for any $\eps\in(0,b]$,
\[
\frac{V(\eps)-V(0)}{\eps}\le -e^{-\rho V(0)-1-\eps^2}.
\]

Letting $\eps\to 0$ we get
\[
\limsup_{\eps\to 0}\frac{V(\eps)-V(0)}{\eps}\le -e^{-\rho V(0)-1},
\]
while (\ref{DVlower}) gives
\[
\frac{V(\eps)-V(0)}{\eps}\ge -\big(1+q(\eps)\big)e^{q(\eps)-1-\rho V(\eps)}.
\]

Letting $\eps\to 0$ and noting that $q(\eps)\to 0$, and $\ V(\eps)\to V(0)$, we have
\[
\liminf_{\eps\to 0}\frac{V(\eps)-V(0)}{\eps}\le -e^{-\rho V(0)-1},
\]
which proves the claim. \hfill$\Box$ \\[2mm]

We conclude observing that since $V\in C\left([0,\infty)\right)$, the general dynamic programming principle is also satisfied. For a proof we refer to \cite{Tou}.

\section{Viscosity solutions and the Hamilton--Jacobi--Bellman equation}

Since the Hamiltonian (\ref{HJB-unbounded}) can take infinite values  we have a singular stochastic control problem and the welfare  function (\ref{ihvf}) should satisfy the proper variational inequality;  see \cite{FS},
Section VIII and \cite{pham}, Section 4. The proof of the  next Lemma, except for the treatment of the boundary conditions, follows the lines of Proposition 4.3.2 of  \cite{pham}. 
\begin{lemma}\label{singular_control_problem}
If $\sigma^2 < \rho +2b$, the welfare function $V$ defined by (\ref{ihvf}) is a continuous constrained viscosity solution of 
\begin{equation}\label{HJB-var}
  \min \Big[ \rho V -\sup_{u\in \mathbb{R}^+ } \left( \dfrac{1}{2}\sigma^{2}x^{2}V_{xx}+( u-bx+\dfrac{x^{2}}{x^{2}+1})V_x +\ln u-x^{2}\right), -V_x   \Big] =0, \quad \text{ in }  [0,\infty).
 \end{equation}
\end{lemma}
\noindent
{\em Proof:} That $V$ is a viscosity solution in $(0, \infty)$ follows as in \cite{pham}, so we omit the details. 
\vskip.05in

Here we briefly discuss the subsolution property at $x=0$. Let $\phi$ be a test function 
such that $V-\phi$ has a maximum at $x=0$ with $V(0)-\phi(0)=0$, and, proceeding by contradiction,  we assume that 
\begin{align}\label{subat0}
  \rho \phi(0) -\sup_{u\in \mathbb{R}^+ } G\big(0,u,\phi'(0),\phi''(0)\big) >0 \,\, \text{ and }\,\, -\phi'(0) >0. 
\end{align}

Since  $-\phi'(0)>0$, the supremum in the above inequality is finite, the Hamiltonian takes the standard form, and the first inequality in \eqref{subat0} becomes
\[
\rho \phi(0)+1+\ln\big(-\phi'(0)\big)>0.
\]

On the other hand, in view of  Proposition \ref{vprop} (iii),  we have \(\rho V(0)+1+\ln\big(-V'(0)\big)=0\), hence \(V'(0)>\phi'(0)\), contradicting that \(V-\phi\) has a maximum at $x=0$.
\vskip.05in

We have now obtained all the necessary material for the proof of Theorem \ref{constrained_vs}.
\vskip.05in

\noindent
{\em  Proof:}  
The fact that (\ref{ihvf}) is a continuous constrained viscosity solution of the equation ${\rm (OHJB)}$  is a direct consequence of the above Lemma  and the fact that inequality  (\ref{DVbounds})  
implies that $p\leq-C$ for any $p\in D^\pm V(x)$, with $x\in (0, \infty)$. The regularity of $V$ in $(0,\infty)$ follows from the classical results for uniformly elliptic 
operators. \hfill$\Box$ \\

Due to the extra regularity of   the welfare function,   the following verification equation holds in $(0, \infty)$, for any optimal pair $(\overline{u}(\cdot)),\overline{x}(\cdot))$, 
\begin{align*}
\rho V(\overline{x}(t))=& \sup_{u\in U } G(\overline{x}(t),u,V_x(\overline{x}(t)),V_{xx}(\overline{x}(t))) \\
=&\left( \dfrac{\overline{x}^{2}(t)}{\overline{x}^{2}(t)+1}-b\overline{x}(t)\right)V_x(\overline{x}(t))-\Big( \ln(-V_x(\overline{x}(t)))+\overline{x}^{2}(t)+1\Big) \\
& +\dfrac{1}{2}\sigma^{2}\overline{x}^{2}(t)V_{xx}(\overline{x}(t)), 
\ \ \ \ \  \ t\in[0,\infty)- a. e.\ ,\ \ \mathbb{P}- a.e.;
\end{align*}
\noindent
 see   \cite{FS} and \cite{JYXZ}.
\vskip.05in

Next we prove the proper  comparison principle for ${\rm (OHJB)}$. The proof is along the lines of  the strategy in \cite{HI}, 
where given a subsolution $u$ and a supersolution $v$ 
of (OHJB),  $u-v$ is a subsolution of the corresponding linearized equation. Then, one concludes by comparing $u-v$ with the appropriate supersolution of the linearized equation; see also 
\cite{DaLL} and \cite{Za}. The difference with the existing results is that, 
due to the presence of the logarithmic term,
the commonly used  functions of simple polynomials do not yield a supersolution of the equation.  
\vskip.05in

Having in mind that we are looking for a viscosity solution that is strictly decreasing 
and satisfies (\ref{DVbounds}),  we prove the following lemma. 

\noindent
\begin{lemma} \label{intermed} Suppose  $u$, $v$ satisfy the assumptions of Theorem \ref{comparison}.
Then $\psi=u-v$ is a subsolution of 
\begin{equation}\label{u-v}
 \rho \psi +bxD\psi-\left(1+c^*\right)|D\psi| -\frac{1}{2}\sigma^2x^2 D^2 \psi = 0 \,\, \mbox{ in } \,\, [0, \infty).
\end{equation}
\end{lemma} 
\noindent
{\em Proof:}
Let  $\bar{x}\ge 0$ a maximum point of  $\psi -\phi$ for some smooth function $\phi$ and set, following \cite{So},
$$\theta(x,y)=\phi(x)+\frac{(x-y+\eps L)^2}{\eps}+\delta(x-\bar{x})^4,$$
where $L,\delta$ are positive constants.

The assumptions on $u,v$ imply that the function $(x,y) \mapsto u(x)-v(y)-\theta(x,y)$ is bounded from above and achieves its maximum at, say, $(x_{\eps},y_{\eps})$. It follows that  $x \mapsto u(x)-v(y_{\eps})-\theta(x,y_{\eps})$ has a local maximum at
$x_{\eps}$ and  $y \mapsto v(y)-u(x_{\eps})+\theta(x_{\eps},y)$ has a local minimum at $y_{\eps}$. Moreover, (see Proposition 3.7 in \cite{CIL}), as  $\eps \rightarrow 0$, 
\begin{equation}\label{doubling}
 \frac{|x_{\eps}-y_{\eps}|^2}{\eps} \rightarrow 0 ,\, x_\eps\rightarrow \bar{x}, \, \mbox{ and }
u(x_{\eps})-v(y_{\eps}) \rightarrow \psi(\bar{x}).
\end{equation}

The inequalities
\[
u(x_\eps)-v(y_\eps)-\theta(x_\eps,y_\eps)\le \psi(\bar{x})-\phi(\bar{x})+v(x_\eps)-v(y_\eps)-\frac{|x_{\eps}-y_{\eps}+\eps L|^2}{\eps}-\delta(x_\eps-\bar{x})^4
\]
and 
\[
u(x_\eps)-v(y_\eps)-\theta(x_\eps,y_\eps)\ge u(\bar{x})-v(\bar{x}+\eps L)-\theta(\bar{x},\bar{x}+\eps L)\ge \psi(\bar{x})-\phi(\bar{x})
\]
together imply that 
\[
\frac{|x_{\eps}-y_{\eps}+\eps L|^2}{\eps}+\delta(x_\eps-\bar{x})^4\le v(x_\eps)-v(y_\eps).
\]

Since $v$ is decreasing we must have $y_\eps>x_\eps.$ In particular, $y_\eps\in(0,\infty)$.

Therefore, setting $p_{\eps}=2\frac{x_{\eps}-y_{\eps}+\eps L}{\eps}$  and $q_\eps= \phi_x(x_{\eps})+4\delta(x_\eps-\bar{x})^3$, Theorem 3.2 in \cite{CIL} implies that,
for every $\alpha >0$, there exist $X, \, Y \in \mathbb{R}$ such that
\begin{equation}\label{sub-super}
 \rho u(x_{\eps})-H(x_{\eps},p_\eps+q_{\eps},X) \leq 0 \mbox{   \,\,and\,\,  } \rho v(y_{\eps})- H(y_{\eps},p_{\eps},Y) \geq 0
\end{equation}
and 
\begin{equation}\label{jet-inequality}
  -(\frac{1}{\alpha}+\|M\|) I \leq \left(
\begin{array}{cc}
X & 0 \\
0 & -Y  
\end{array} \right) \leq M+\alpha M^2
\end{equation}
with $M=D^2\theta(x_{\eps},y_{\eps}).$

By subtracting the two inequalities in  (\ref{sub-super}) we obtain 
\begin{align}\label{comp-inequ_1}
 \rho & u(x_{\eps})-\rho v(y_{\eps}) +b(x_{\eps}-y_{\eps})p_{\eps}+\Big(\frac{y_{\eps}^2}{y_{\eps}^2+1}-\frac{x_{\eps}^2}{x_{\eps}^2+1}\Big)p_{\eps}
 - \Big(\frac{x_{\eps}^2}{x_{\eps}^2+1}-bx_{\eps} \Big) q_\eps  \\
 & +\ln{(-p_{\eps}- q_\eps)}- \ln{(-p_\eps)} \nonumber + x_{\eps}^2- y_{\eps}^2  -\frac{1}{2}\sigma^2x_{\eps}^2 X + \frac{1}{2}\sigma^2y_{\eps}^2 Y \leq 0. \nonumber
\end{align}

Our assumption on $u$ implies that $p_\eps+q_\eps\le -1/c^*$. Thus, the difference of the logarithmic terms in the above inequality can be estimated from below as
\[
\ln\Big(\frac{p_\eps+q_\eps}{p_\eps}\Big)\ge \frac{q_\eps}{p_\eps+q_\eps}\ge -c^*|q_\eps|,
\]
and inequality (\ref{comp-inequ_1}) gives
\begin{align}\label{comp-inequ_2}
\frac{1}{2}\sigma^2x_{\eps}^2 X + \frac{1}{2}\sigma^2y_{\eps}^2 Y&\ge\rho u(x_{\eps})-\rho v(y_{\eps}) + bx_\eps q_\eps- (1+c^*)|q_\eps| +\nonumber\\
 &\qquad p_\eps(x_\eps-y_\eps)\Big(b-\frac{x_\eps+y_\eps}{(1+x_\eps^2)(1+y_\eps^2)}\Big)+x_\eps^2-y_\eps^2
\end{align}

On the other hand, the right-hand-side in (\ref{jet-inequality}) yields
\begin{equation}\label{X-Y}
 \frac{1}{2}\sigma^2x_{\eps}^2 X - \frac{1}{2}\sigma^2y_{\eps}^2 Y \leq \frac{1}{2}\sigma^2x_{\eps}^2 \big(\phi_{xx}(x_{\eps})+12\delta(x_\eps-\bar{x})^2\big)
 +\frac{\sigma^2}{\eps}(x_{\eps}-y_{\eps})^2+m(\frac{\alpha}{\eps^2}), 
\end{equation}
with $m$ a modulus of continuity independent of $\alpha$, $\eps$. 
\vskip.05in

By combining \eqref{comp-inequ_2} with \eqref{X-Y}, we conclude the proof taking first $\alpha  \rightarrow 0$, then $\eps \rightarrow 0$ and using \eqref{doubling}. \hfill$\Box$ \\

\noindent
 We continue with the \\[1mm]
 \noindent
{\it Proof of Theorem \ref{comparison}:} The main step is the construction of a solution of the linearized equation. For this, we  consider the ode 
\begin{equation}\label{ode1}
\rho w +\big(bx-(1+c^*)\big)w'-\frac{1}{2}\sigma^2x^2w''=0,
\end{equation}
which has a solution of the form 
\begin{equation}\label{sol}
 w(x)=x^{-k}\mathcal{J}(\frac{2+2c^*}{\sigma^2x}),
\end{equation}
where $k$ is a root of  
\begin{equation}\label{root}
 k^2+\Big(1+\frac{2b}{\sigma^2}\Big)k-\frac{2\rho}{\sigma^2}=0
\end{equation}
and   $\mathcal{J}$ a solution of the degenerate hypergeometric equation 
\begin{equation}\label{confluent}
xy''+(\tilde{b}-x)y'-\tilde{a}y=0
\end{equation}
with $\tilde{a}=k$ and $\tilde{b}=2(k+1+b/\sigma^2)$.

Since we are looking for a solution of \eqref{ode1} with superquadratic growth at $+\infty$, we choose $k$ to be the negative root of \eqref{root}. The assumption $\sigma^2 < \rho+2b$ implies $-k>2$.

We further choose $\mathcal{J}$ to be the Tricomi solution of \eqref{confluent} which satisfies
\[
\mathcal{J}(0)>0\qquad\text{and} \qquad \mathcal{J}(x)=x^{-k}\big(1+\frac{2\rho}{\sigma^2 x}+o(x^{-1})\big)\quad\text{as }x\to\infty.
\]

With this choice, the function $w$ defined in \eqref{sol} for $x>0$ and by continuity at $x=0$, satisfies $w(0), w'(0) >0$ and $w(x)\sim \mathcal{J}(0) x^{-k}$, as $x\to\infty$. 

Note that $w$ is increasing in $[0,\infty)$ since it would otherwise have a positive local maximum and this is impossible by \eqref{ode1}. In particular, $w$ satisfies \eqref{u-v}.

Set now $\psi=u-v$ and consider $\epsilon>0$. Since $\psi-\epsilon w<0$  in a neighborhood of infinity, there exists $x^\epsilon \in [0,\infty)$ such that
\[
\max_{x\ge 0} \big(\psi(x)-\epsilon w(x)\big)=\psi(x^\epsilon)-\epsilon w(x^\epsilon).
\]

By Lemma \ref{intermed} $\psi $ is a subsolution of \eqref{u-v}. We now use $\epsilon w$ as a test function to find that
\begin{equation*}
0\ge \rho \psi({x}^\epsilon) +\epsilon b{x}^\epsilon w({x}^\epsilon)-\epsilon\left(1+c^*\right)|w'({x}^\epsilon))| -\frac{1}{2}\sigma^2{(x^\epsilon)}^2 w''({x}^\epsilon) = \psi({x}^\epsilon)-\epsilon w(x^\epsilon).
\end{equation*}

Hence, $\psi(x)\le \epsilon w(x)$ for all $x\in [0,\infty)$. Since $\epsilon$ is arbitrary, this proves the claim. \hfill$\Box$ \\[1mm]

The stability property of viscosity solutions yields the following theorem. 

\begin{thm}\label{constrained_vs_det}  As  $\sigma \rightarrow 0$, the welfare function $V$ defined by (\ref{ihvf}) converges locally uniformly to the constrained viscosity solution $V^{(d)}$   of the deterministic 
shallow lake equation in $[0, \infty)$,
$$\rho V^{(d)}=\left( \dfrac{x^{2}}{x^{2}+1}-bx\right)V_{x}^{(d)}-\Big( \ln(-V_{x}^{(d)})+x^{2}+1\Big). \leqno\mathrm{(OHJB_d)}$$
\end{thm}
\vskip.05in
We next prove Theorem  \ref{asymptotic} that describes the asymptotic behaviour of $V$ as $x \rightarrow \infty$. The proof is based on a scaling argument and the stability properties of the viscosity solutions.\\

\noindent
{\it Proof of Theorem \ref{asymptotic}: } We write  $V$ as 
\begin{equation*}
 V(x)= -A\left(x+\frac{1}{b+\rho}\right)^2-\frac{1}{\rho}\ln \left(2A(x+\frac{1}{b+\rho})\right)+K+v(x). 
\end{equation*}

Straightforwad calculations yield that   $v$ is a viscosity solution in $(0, \infty)$ of the equation
\begin{multline}
\rho v+\left(bx-\dfrac{x^{2}}{x^{2}+1}\right)v'+\ln\left(1+\frac{1-\rho\big(x+\frac{1}{b+\rho}\big)v'}{2A\rho\left(x+\frac{1}{b+\rho}\right)^2} \right)-\dfrac{1}{2}\sigma^{2}x^{2}v''+f=0,
\end{multline}
where 
\[
f(x)=\frac{b+\frac{\sigma^2}{2}+(b+\rho)\frac{x^2}{1+x^2}}{\rho\big(1+x(b+\rho)\big)}+\frac{\sigma^2x(b+\rho)}{2\rho\big(1+x(b+\rho)\big)^2}-2A\frac{\big(1+x(b+\rho)\big)}{1+x^2}.
\]

Note $f$ is smooth on $[0,\infty)$ and vanishes as $x\to\infty$.
\vskip.05in

 Let $v_\lambda(y)=v(\frac{y}{\lambda})$ and observe that, if $v_\lambda(1) \rightarrow 0$ as $\lambda \rightarrow 0$, then $v(x) \rightarrow 0$ as $x 
 \rightarrow \infty$. It turns out that  $v_\lambda$ solves
 \begin{multline*}
\rho v_\la+ \left(bx-\dfrac{\la x^{2}}{x^{2}+\la^2}\right)v'_\la+\ln\left(1+\frac{\la^2\big(1-\rho\big(x+\frac{\la}{b+\rho}\big)v'_\la\big)}{2A\rho\left(x+\frac{\la}{b+\rho}\right)^2} \right)-\dfrac{1}{2}\sigma^{2}x^{2}v''_\la+f\big(\frac{x}{\la}\big)=0.
 \end{multline*}

Since, by (\ref{Vbounds}) $v_\la $ is uniformly bounded, we consider the half-relaxed limits $v^*(y)=\limsup_{x \rightarrow y, \la\rightarrow 0}v_\la(x)$ and $v_*(y)=\liminf_{x \rightarrow y, \la\rightarrow 0}v_\la(x)$ 
 in $(0, \infty)$, which are (see \cite{BP}) respectively sub- and super-solutions of 
\begin{equation}
 \rho w +byw'-\frac{1}{2} \sigma^2 y^2 w''=0.
\end{equation}

It is easy to check that for any $y>0$ we have $v^*(y)=\limsup_{x\to\infty}v(x)$ and $v_*(y)=\liminf_{x\to\infty}v(x)$. 
\vskip.05in
The subsolution property of $v^*$ and the supersolution property of $v_*$ give 
\[
\limsup_{x\to\infty}v(x)\le 0\le\liminf_{x\to\infty}v(x)\le \limsup_{x\to\infty}v(x).
\]\qed

\section{A numerical scheme and optimal dynamics}

 A general argument to prove the convergence of monotone schemes for viscosity solutions of fully nonlinear second-order elliptic or parabolic, possibly degenerate, partial differential equations has been introduced in \cite{BS}. 
 Their methodology  has been extensively
 used to approximate solutions to first-order equations, see for example  \cite{RT},  \cite{Se2}, \cite{KMS}, \cite{Qi}. 
 \vskip.05in
 
 On the other hand, it is not always possible to construct  monotone schemes for second-order equations in their full generality. However, various types of nonlinear second-order equations have been approximated
 via monotone schemes based on \cite{BS}; see, for example,   \cite{OF}, \cite{BJ}, \cite{BZ}, \cite{FO}.
 \vskip.05in
 
 Next, following \cite{KosZoh} which considered the deterministic problem, we construct a monotone finite difference scheme to approximate numerically the welfare function and recover numerically the stochastic optimal dynamics. 
\vskip.05in  

Let
$\Delta x$ denote the step size of a uniform partition $0 = x_0 < x_1 < \ldots < x_{N-1} < x_N = l$ of $[0, l]$ for 
$l > 0$ sufficiently large. Having in mind (\ref{DVbounds}),  if $V_{i}$ is the approximation of $V$ at $x_i$,  we employ a backward finite difference discretization  to approximate the first derivative in the linear term of the (OHJB),  
a forward finite difference discretization for the derivative in the logarithmic term and a central finite difference scheme to approximate the second derivative. 

\vskip.05in
These considerations yield, for $i = 1,\ldots, N-1$,   the approximate equation
\begin{multline}\label{DHJB}
 V_i - \frac{1}{\rho} \Big(\frac{x_i^2}{x_i^2+1} - bx_i \Big) \frac{ V_{i} - V_{i-1} }{ \Delta x} + \frac{1}{\rho}
\left[x_i^2 + 1 +\ln \left( - \frac{V_{i+1} - V_i}{\Delta x}\right)\right] \\
- \frac{\sigma^2}{2\rho}
\frac{ V_{i+1} + V_{i-1} - 2 V_i}{(\Delta x)^2}=0.
\end{multline}

Setting 
\begin{multline}
 g(x,w,c,d)=\Big[(\Delta x)^2 - \frac{1}{\rho}\Big(\frac{x^2}{x^2+1} - bx \Big) \Delta x +\frac{\sigma^2}{\rho}\Big]w+ \\ 
 \frac{1}{\rho}
(x^2 + 1 ) (\Delta x)^2+\frac{1}{\rho}(\Delta x)^2\ln \left( - \frac{c - w}{\Delta x}\right)+\frac{1}{\rho}\Delta x \Big(\frac{x^2}{x^2+1} - bx \Big) d- \frac{\sigma^2}{2\rho}
( c + d ),
\end{multline}
 the numerical  approximation of $V$ satisfies 
\begin{equation}\label{fdscheme}
 g(x_i,V_i,V_{i+1},V_{i-1})=0, \mbox{\,\,\,for\,\, \, i = 1,\ldots, N-1 },
\end{equation}
and the consistency is immediate.
\vskip.05in

For the monotonicity we observe that for  two different approximation grid vectors $(U_0, \ldots, U_N)$ and $(V_0, \ldots, V_N)$ with $U_i \geq V_i$ and $U_i=V_i=w$, we have 
\begin{equation}
 g(x_i,w,U_{i+1},U_{i-1})\leq g(x_i,w,V_{i+1},V_{i-1}),
\end{equation}
provided $\Delta x (\frac{x^2}{x^2+1} - bx)\leq \sigma^2/2$. This condition is satisfied if $b \geq 0.5$ or if we take $\Delta x\leq\sigma^2/2$.
\vskip.05in
Since the welfare function solves a state constraint problem the equation is satisfied on the left boundary point.
\vskip.05in

It follows that  the  numerical scheme is monotone in the sense of \cite{BS} and converges to the unique constrained
viscosity solution. 
\vskip.05in

 Since the computational  domain of the problem is finite, a boundary condition has to be imposed at $x=l$, for $l$ sufficiently large, by exploiting  the asymptotic behaviour of the welfare function $V$ as $x\to+\infty$. The boundary condition at
the right endpoint $x_N$ is provided by  the asymptotic estimate \eqref{asympto_behav}.  
\vskip.05in

The scheme above suggests a numerical algorithm for the computation of optimal dynamics governing the shallow lake problem. In this direction, the nondegeneracy of the shallow lake equation in $(0,\infty)$ induces  extra regularity for the function $V$ in $(0,\infty)$. Hence, the optimal dynamics for the shallow lake problem  are described by 
 \begin{equation}\label{opt_dyn}
 \left\{ \begin{array}{l} d\bar{x}(t)=\left(-\dfrac{1}{V'(\bar{x}(t))} -b \bar{x}(t)+\dfrac{\bar{x}^{2}(t)}{\bar{x}^{2}(t)+1}\right)dt+\sigma \bar{x}(t)dW(t), \\
   \bar{x}(0)=x       
  \end{array} \right.
\end{equation}
 
Using the numerical representation  of $V$ via \eqref{fdscheme} and  properly discretizing the SDE (\ref{opt_dyn}), we can  reconstruct numerically the optimal dynamics. 
This is a direct approach to investigate numerically the stochastic properties of the optimal dynamics of the shallow lake problem for the various parameters $\rho$, $b$, $c$, $\sigma$ of the problem.
\vskip.05in

The exact numerical algorithm for the computation of the constrained viscosity solution along with the numerical study of the optimal dynamics and their stochastic properties 
 for various $\sigma$'s will be presented elsewhere.

\appendix
\section{}

%
\begin{lemma}\label{Mest}
Assume that $f$ is a positive $\PP$-a.s. locally integrable $\mathcal{F}_t$ and let $M_t(f)$ be defined as in (\ref{Mf}). 
Then,\\[2mm]
(i) \quad $\displaystyle \EE\left[\int_0^\infty e^{-\rho t}M_t(f)\ dt\right]=\frac{1}{\rho+b}\ \EE\left[\int_0^\infty e^{-\rho t}f(t)\ dt\right].$\\[2mm]
(ii) \quad $\displaystyle
\EE\left[\int_0^\infty e^{-\rho t} Z_t M_t(f)\ dt\right]=\begin{cases}A\ \EE\left[\int_0^\infty e^{-\rho t} Z_t f(t)\ dt\right] &\text{ if }  \sigma^2<\rho+2b,\\
										 \infty & \text{ if } \sigma^2\ge\rho+2b. \end{cases}
$\\[2mm]
(iii) \quad$\displaystyle
\EE\left[\int_0^\infty e^{-\rho t} M_t^2(f)\ dt\right]=\begin{cases}2A\ \EE\left[\int_0^\infty e^{-\rho t} f(t)M_t(f)\ dt\right] &\text{ if }  \sigma^2<\rho+2b,\\
										 \infty & \text{ if } \sigma^2\ge\rho+2b. \end{cases}			
$
\end{lemma}
\noindent
{\em Proof:} (i) Since $f$ is $\mathcal{F}_t$-adapted we have
\[
\EE\big[M_t(f)\big]=\EE\left[\int_0^t \EE\big[Z_t\big| \mathcal{F}_s\big] \frac{f(s)}{Z_s}\ ds\right]=\int_0^t e^{-b(t-s)}\EE\big[f(s)\big]\ ds.
\]

Therefore,
\begin{align*}
\EE\left[\int_0^\infty e^{-\rho t}M_t(f)\ dt\right]&=\int_0^\infty e^{bs}\EE\big[f(s)\big]\int_s^\infty e^{-(\rho+b)t}dt\ ds\\&=\frac{1}{\rho+b}\ \EE\left[\int_0^\infty e^{-\rho t}f(t)\ dt\right].
\end{align*}
\noindent

(ii) Conditioning first on $\mathcal{F}_s$ we have
\begin{align*}
\EE\big[Z_tM_t(f)\big] &= \EE\left[\int_0^t\EE\big[Z_t^2\big|\mathcal{F}_s\big] \frac{f(s)}{Z_s}\  ds\right]=\int_0^t e^{(\sigma^2-2b)(t-s)}\EE\left[Z_s{f(s)}\right] \  ds,
\end{align*}
and, hence, 
\begin{align*}
\EE\left[\int_0^\infty e^{-\rho t} Z_t M_t(f)\ dt\right]=&\int_0^\infty e^{(\sigma^2-2b)s} \EE\left[Z_sf(s)\right]\int_s^\infty e^{-(\rho+2b-\sigma^2)t}\ dt \ ds\\
=&\begin{cases}A\ \EE\left[\int_0^\infty e^{-\rho t} Z_t f(t)\ dt\right] &\text{ if }  \sigma^2<\rho+2b\\										 \infty & \text{ if } \sigma^2\ge\rho+2b. \end{cases}
\end{align*}

(iii) By Fubini's theorem we have
\begin{align*}
\EE\big[M_t^2(f)\big]&=2\ \EE\left[\int_0^t\int_s^t \frac{Z_t^2}{Z_sZ_q}f(s)f(q)\ dq\ ds\right]\\
&=2\ \EE\left[\int_0^t\int_s^t \EE\big[Z_t^2\big|\mathcal{F}_q\big]\frac{1}{Z_sZ_q}f(s)f(q)\ dq\ ds\right]\\
&=2 \int_0^t\int_s^t e^{(\sigma^2-2b)(t-q)}\EE\left[\frac{Z_q}{Z_s}f(s)f(q)\right]\ dq\ ds\\
&=2 \int_0^t e^{(\sigma^2-2b)(t-q)}\EE\left[f(q)M_q(f)\right]\ dq,
\end{align*}
and, therefore,
\begin{align*}
\EE\left[\int_0^\infty\hspace{-2mm} e^{-\rho t}M_t^2(f)\ dt\right]&=\int_0^\infty \hspace{-2mm} e^{(2b-\sigma^2)q}\EE\big[f(q)M_q(f)\big]\int_q^\infty \hspace{-2mm} e^{-(\rho+2b-\sigma^2)t}dt\ dq\\
&=\begin{cases}2A\ \EE\left[\int_0^\infty e^{-\rho t} f(t)M_t(f)\ dt\right] &\text{ if }  \sigma^2<\rho+2b,\\									 \infty & \text{ if } \sigma^2\ge\rho+2b.\ \Box \end{cases}
\end{align*}

\end{document}